\documentclass[a4paper,11pt]{article}
\usepackage[T1]{fontenc}
\usepackage[utf8]{inputenc}
\DeclareUnicodeCharacter{03BB}{\ensuremath{\lambda}}
\usepackage[bitstream-charter]{mathdesign}

\usepackage{authblk}
\usepackage{amsmath}
\usepackage{mathtools}
\usepackage{cite}
\usepackage{enumitem}
\usepackage{xcolor}
\usepackage{url}
\usepackage[all]{nowidow}
\usepackage{float}
\usepackage{graphicx}
\usepackage{fancyvrb}
\usepackage{vmargin}
\setmarginsrb{2cm}{2cm}{2cm}{2cm}{0mm}{0mm}{0mm}{10mm}
\usepackage{hyperref}
\hypersetup{colorlinks=true,allcolors=blue}

\title{\vspace*{-1.5cm} \bfseries Matrix diagonalization and\\ singular value decomposition:\\ Static SageMath and dynamic ChatGPT juxtaposed}
\author{Natanael Karjanto}
\affil{Department of Mathematics, University College, Natural Science Campus\par Sungkyunkwan University, Suwon~16419, Republic of Korea}

\date{\vspace*{-0.5cm} \scriptsize Updated \today}

\begin{document}
\maketitle

\begin{abstract}
We investigated some difficulties that students often face when studying linear algebra at the undergraduate level, and identified some common mistakes and difficulties they often encountered when dealing with topics that require algorithmic thinking skills such as matrix factorization. In particular, we focused on (orthogonal) diagonalization and singular value decomposition (SVD). We also offered the possibility of exploring these topics using SageMath, a Python-based free open software computer algebra system (CAS) that has been identified to be useful for assisting many students in the computational process even though its output is static by nature. We then explored dynamic ChatGPT by inquiring the chatbot about the topic, either by asking to provide an example or to solve a problem, that is by constructing an (orthogonal) diagonalization or SVD from a particular matrix. By consolidating essential concepts in linear algebra and improving computational skills through effective practice, mastering these topics would become easier and mistakes could be minimized. Static SageMath, in particular, is a great aid for calculation confirmation and handling tedious computations. Although dynamic ChatGPT is relatively unreliable for solving problems in linear algebra, the mistakes it produces could become a valuable tool for improving critical thinking skills.
\end{abstract}

\section{Introduction}

Linear algebra is a branch of mathematics that deals with linear equations, linear functions, and their representations in vectors and matrices. It involves the study of vector spaces, linear transformations, matrices, determinants, eigenvalues, and eigenvectors. Linear algebra is a fundamental tool in many areas of mathematics, including geometry, calculus, optimization, and numerical analysis. It also has practical applications in various fields such as physics, engineering, computer science, and economics~\cite{strang2006linear,friedberg2018linear,axler2015linear,anton2019elementary,larson2016elementary,andrilli2022elementary,hill2014elementary,kolman2017elementary,bretshcer2018linear,kuttler2021elementary,baker2014linear,spence2017elementary,lay2021linear,mirsky2013an,johnson2017introduction,hoffman2015linear,lang1997introduction,blyth2002,boyd2018introduction,lax2007linear}.

Similar to other fields of mathematics, algorithmic thinking in linear algebra is a crucial skill in understanding and mastering the materials. Often also called computational thinking, algorithmic thinking is the process of breaking down complex mathematical problems into smaller, more manageable parts, and then solving them systematically using a sequence of steps or algorithms. This process involves identifying patterns, constructing algorithms, and developing logical and analytical skills in solving these problems~\cite{barcelos2018mathematics,stephens2020computational,sung2017introducting,weintrop2016defining}. Teaching and learning linear algebra for both mathematics and non-mathematics majors provides an opportunity for cultivating and embracing problem-solving, logical reasoning, and computational thinking skills, the very skills that are essential in various areas of science and engineering requiring a solid knowledge of linear algebra, such as operational research, computer science, data science, and machine learning, among others~\cite{axler2015linear,strang2006linear,boyd2004convex,trefethen1997numerical,bennett2021linear,cohen2022practical,golub2013matrix,klein2013coding,strang2019linear,aggarwal2020linear,yoshida2021linear}. Strengthening computational thinking among future teachers will also be useful when they eventually train the future generations of students~\cite{gadanidis2017computational,yadav2017computational}. 

In addition to being able to manipulate a matrix algebraically and to perform matrix operations, some essential concepts in linear algebra include but are not limited to understanding vector spaces and subspaces, applying linear transformation, working with norms and inner products, finding least-squares solutions, and the ability to perform matrix factorization. In this study, we focused on the latter. Sometimes also called decomposition, factorization in linear algebra refers to the process of decomposing a matrix into several, often simpler, matrices that can be analyzed and manipulated more easily. These simpler matrices may have special properties or structures that make them easier to analyze or compute with. Mastering algorithmic thinking in performing matrix factorization is essential for every linear algebra's learner because this topic not only requires other basic concepts in linear algebra but is also an essential building block for understanding other more complicated topics in linear algebra.

There are various types of matrix decomposition, such as LU decomposition, QR decomposition, diagonalization, orthogonal diagonalization, singular value decomposition (SVD), and eigendecomposition. The latter refers to the factorization of a matrix into a canonical form, where the matrix is represented by its eigenvalues and eigenvectors. Each type of decomposition has its own specific properties and applications. For example, LU decomposition is useful for solving systems of linear equations, QR decomposition is used for least-squares problems, SVD is used for data compression and feature extraction, and eigendecomposition is used for analyzing the behavior of linear operators. Whereas matrix diagonalization only applies to square matrices and orthogonal diagonalization only applies to a special class of symmetric square matrices, SVD can be implemented for any size of the rectangular matrix, which does not necessarily have to be square.

In short, the SVD is a factorization method in linear algebra that decomposes a real- or complex-valued matrix into three components: a diagonal matrix of singular values and two unitary (orthogonal) matrices. Mathematical applications of SVD encompass the calculation of the (Moore-Penrose) pseudoinverse, approximating matrix, and determining the rank, range, and null space of a matrix. SVD is a widely utilized technique in data analysis and machine learning for reducing the dimensionality of data and extracting important features. SVD is also used in a variety of applications in science and engineering, including signal and image processing, text mining, data least-squares fitting, process control, and recommendation systems~\cite{abdi2007singular,akritas2004applications,andrews1976singular,biglieri1989some,bisgard2021analysis,freire1988application,henry1992singular,klema1980the,lange2010singular,maciejewski1989the,mandel1982use}. The singular values obtained from SVD represent the importance of each feature in the original data and can be used to reconstruct the data with different levels of accuracy.

A common curriculum in linear algebra outlines that students ought to know how to calculate both matrix diagonalization and orthogonal diagonalization before introducing them to SVD. However, we observed that students often encountered many difficulties with these decompositions. For orthogonal diagonalization, several students often forget to normalize the associated eigenvectors to become unit vectors. In other cases, when one eigenvalue yields two associated eigenvectors, some students also forget to transform the set of eigenvectors into an orthonormal set, which can be easily done by orthogonal projection via the Gram-Schmidt process, keeping one of them while projecting the other.

For SVD, a common mistake is in finding the singular values of a matrix. Let $A$ be a matrix, and instead of finding the eigenalues of $A^T A$, some students calculated the eigenvalues of $A$, took their square root, and designated them as singular values. Another difficulty occurs when the matrix does not have a full rank, for which the students need to find the missing one or more orthonormal eigenvectors in one of the matrices by implementing the orthogonal property or cross product operation. In the absence of these eigenvectors, the resulting factorization is called ``reduced SVD'' instead of (full) SVD. 

To the best of our knowledge, this issue of common mistakes and learners' struggles in understanding (orthogonal) diagonalization and SVD has not been fully addressed in the body of published literature. The closest article that is tangentially related to our study is Yildiz Ulus' (2013) paper, in which the author investigated teaching diagonalization using advanced calculators and observed that such technological tools are beneficial for learners' acquisition of algorithmic mathematical knowledge for this particular topic of linear algebra~\cite{yildiz2013teaching}. Although Lazar (2012) argued in his master's thesis that a solid understanding of fundamental concepts in linear algebra is essential for mastering SVD, the author did not inquire the participants in his study regarding the difficulty or common mistakes they encountered when studying the topic~\cite{lazar2012teaching}. 

This article attempts to fill this gap by investigating common mistakes and challenges when students learn linear algebra, particularly (orthogonal) diagonalization and SVD. Furthermore, because teaching and learning in contemporary mathematics cannot be dismantled by utilizing technological tools, we are also interested in investigating whether computer software or another newly arrived technology can assist and enhance students' learning, instead of disrupting and marring it. Our scope for the computer algebra system (CAS) is SageMath, which has a static nature, and we will consider ChatGPT as an artificial intelligence (AI)-assisted tool, which has dynamic interaction features. In the following, we provide a brief overview of static CAS SageMath and dynamic ChatGPT.
 
What is a CAS? A CAS is a software program that allows the manipulation and computation of mathematical expressions and symbols, including algebraic equations, calculus, and other mathematical functions. It is designed to perform symbolic manipulation, numerical computations, and graphics, as well as to provide tools for solving equations, manipulating matrices and vectors, and performing other mathematical operations. CAS is commonly used in scientific research, engineering, and education, and can be used to solve complex mathematical problems that may be too difficult or time-consuming to solve by hand~\cite{buchberger1983computer,chudnovsky1990computers,cohen2002computer,davenport1993computer,fey2003computer,fuchs2003computer,geddes1992algorithms,grabmeier2003computer,harper1991guide,koepf2021computer,lamagna2019computer,simon1996interoprability,soiffer1992the,von2013modern,wester1999computer}. Examples of popular CAS include (wx)Maxima~\cite{hannan2015wxmaxima,kanagasabapathy2018intro,kadry2019mathematics,senese2019symbolic,timberlake2016classical},
Wolfram~Mathematica~\cite{abell2022mathematica,maeder1991programming,trott2007the,wagon1999mathematica,wolfram2003the}, Maplesoft~Maple~\cite{abell2005maple,fox2019advanced,fox2020advanced,garvan2001the,heck2003introduction}, and Mathworks~MATLAB		~\cite{attaway2022matlab,davis2010matlab,etter2017introduction,higham2017matlab,valentine2022essential}.

What is SageMath? SageMath (also known as Sage) is a free open-source mathematics software system that uses Python as its primary programming language. It aims to provide an alternative to commercial mathematical software systems, such as Wolfram~Mathematica and Maplesoft~Maple, while also providing an interface to other popular mathematics software systems, such as Mathworks~MATLAB and GAP, a system for computational discrete algebra. \emph{SageMath} has a wide range of capabilities, including algebraic and numerical computations, graphics, symbolic manipulation, and combinatorics. Its development is community-driven and released under the GNU General Public License~\cite{bard2015sage,finch2011sage,joyner2008sage,kumar2022basics,szabo2017sage,zimmermann2018computational}.

Interestingly, Sagemath can not only be used for linear algebra but also has very strong support for the subject, with many built-in functions and capabilities for working with matrices, vectors, and linear transformations. SageMath provides a variety of algorithms and methods for solving linear algebra problems, such as performing matrix operations, solving linear systems of equations, finding eigenvalues and eigenvectors, and computing matrix decompositions (e.g., LU, QR, SVD). SageMath is a powerful tool for symbolic linear algebra, which allows for exact computations with variables and expressions. Despite being static in nature, we are convinced that SageMath can be useful for teaching and learning linear algebra, not only for assisting computational tasks but also for helping to check whether our hand calculations were correct. In addition, the time saved for computation can be utilized for other purposes, such as understanding deeper mathematical concepts in linear algebra or exploring various problems in the topic.

What is ChatGPT? ChatGPT (Chat Generative Pre-trained Transformer) is a language model developed by OpenAI, an American AI research laboratory consisting of a team of researchers and engineers dedicated to creating safe and beneficial AI. Although OpenAI was founded in 2015, ChatGPT, as one of its products, was launched as a prototype in November 2022. Amazingly, it reached one million users within five days after its launch. An AI tool built on top of the unsupervised transformer language model GPT-3, ChatGPT was trained on a large dataset of text and could generate responses to questions, write coherent paragraphs, and even carry on a conversation with users. Its purpose is to assist with various tasks, such as answering questions, providing information, or generating text in a conversational manner. 

Interestingly, as an AI language model, ChatGPT has been trained on a large corpus of texts, including mathematical concepts and problems. Thus, it is capable of solving mathematics problems, at least that is what it claims. ChatGPT admits that its ability to solve mathematical problems may be limited by its training data, although it can certainly provide assistance and guidance on various mathematical topics. Specifically, when ask whether it can solve problems in linear algebra, ChatGPT claims of ``possessing knowledge and understanding of linear algebra concepts and can provide solutions to problems in this field.'' During the past few months, the number of published articles related to ChatGPT and its capabilities has increased steadily, including many that appeared in various preprint servers. The following examples cover only a few articles on ChatGPT related to mathematics. 

Recently, Frieder et al. (2023) investigated ChatGPT's mathematical capabilities by asking a wide range of questions, and although it understood the question, it often failed to provide correct solutions. They concluded that the mathematical abilities of ChatGPT were significantly lower than those of average mathematics graduate students~\cite{frieder2023mathematical}. Shakarian et al. (2023) also evaluated ChatGPT for mathematical word problems and discovered that its performance was dramatically altered based on the requirement to show its work. When the chatbot was not required to reveal the complete solution, it failed 84\% of the time, whereas it failed only 20\% when revealing the detailed solution was requested~\cite{shakarian2023an}. Azaria (2022) tested the numerical literacy of ChatGPT and noticed that when it comes to using number digits, ChatGPT is rather biased, with $7$ being the most frequent digit generated by the machine, which also turns out many people's favorite number~\cite{azaria2022chatgpt}. Borji (2023) focused on ChatGPT's failures in mathematics, including arithmetic, logic, and reasoning, among others~\cite{borji2023a}.

In terms of output production, SageMath and other CASs tend to be static, whereas ChatGPT and its competitors are dynamic. The latter could provide a step-by-step explanation of a solution to a particular problem. For matrix factorization in linear algebra, it is essential for learners to understand not only the technical details of a calculation but also to grasp the algorithm behind any particular problem task. Certainly, matrix factorization is a broad topic in itself, and any attempt to discuss other types of matrix factorization should be addressed separately elsewhere. Presently, our focal point for the topic of ``matrix factorization'' only covers diagonalization, orthogonal diagonalization, and SVD. Based on the above literature, we consider the following research questions:
\begin{itemize}[leftmargin=1em]
\item What are some common mistakes and difficulties that students encounter when learning matrix factorization in linear algebra? 

\item How can static CAS, such as SageMath, assist students in learning linear algebra, particularly in (orthogonal) diagonalization and SVD?

\item Can we rely on the dynamic ChatGPT to better understand (orthogonal) diagonalization and SVD?
\end{itemize}

The remainder of this article is organized as follows. After this introduction, Section~\ref{common} features some common mistakes that students often make when learning (orthogonal) diagonalization and SVD. Section~\ref{sagemathsection} continues with SageMath and its ability to assist students in learning (orthogonal) diagonalization and SVD. Section~\ref{chatgptsection} features some examples in which we asked the chatbot to solve problems related to (orthogonal) diagonalization and SVD. We also discuss where ChatGPT makes mistakes and encounters troubles in completing assigned tasks. Finally, Section~\ref{discussandconcludesection} discusses and concludes the study. 

\section{Common mistakes in matrix factorization}	\label{common}

\subsection{Matrix diagonalization}		\label{madi21}

In this section, we consider some common mistakes that students encounter when learning to diagonalize both non-symmetric and symmetric matrices. We also consider similar aspects of SVD. For the former, if $A$ is an $n \times n$ square matrix, then $A$ is diagonalizable whenever it possesses $n$ linearly independent eigenvectors. Usually, the students did not encounter difficulties in diagonalizing a non-symmetric diagonalizable matrix as long as they understood the procedure for finding it. A typical algorithm starts with finding eigenvalues and their associated eigenvectors, constructing a diagonal matrix~$D$, constructing an invertible matrix~$P$, finding its inverse~$P^{-1}$, and expressing the diagonalization, that is, $A = P D P^{-1}$. The following example illustrates this algorithm.

Given a $3 \times 3$ matrix $A$ as follows:
\begin{equation*}
A = \begin{bmatrix*}[r]
1 & 2 & 2 \\
0 &	2 & 1 \\
0 & 1 & 2
\end{bmatrix*}.
\end{equation*}
To diagonalize a matrix, we must determine its eigenvalues and their corresponding eigenvectors. The former can be obtained by solving the characteristic equation $\det \left(A - \lambda I \right) = 0$ or $\left(\lambda - 1\right) \left(\lambda - 3 \right)\left(\lambda - 1 \right) = 0$, which gives $\lambda_1 = 1$, $\lambda_2 = 3$, and $\lambda_3 = 1$. Note that some students might attempt to express the characteristic equation in cubic form, that is, $\lambda^3 - 5\lambda^2 + 7\lambda - 3 = 0$, and solve this for $\lambda$. Although in other cases, this might be a necessary step and even inevitable, we do not have to do this particular step in this specific example. Up to this point, we can form the diagonal matrix $D$ by placing the eigenvalues along the diagonal:
\begin{equation*}
D = \begin{bmatrix*}[r]
1 & 0 & 0 \\
0 &	3 & 0 \\
0 & 0 & 1
\end{bmatrix*}.
\end{equation*}
However, this is not the only option. Depending on the choice of eigenvalues designated as the first, second, third, and so on, we might obtain a different expression for $D$. Coincidentally, both the algebraic multiplicity and geometric multiplicity for the eigenvalue $\lambda = 1$ from the matrix in this example are equal, that is, two. In this case, we are guaranteed to obtain three linearly independent eigenvectors and thus, the matrix is diagonalizable.

Let $\mathbf{p}_1$, $\mathbf{p}_2$, and $\mathbf{p}_3$ be the corresponding eigenvectors, we found that
\begin{equation*}
\mathbf{p}_1 = 
\begin{bmatrix*}[r]
0 \\ 1 \\ -1
\end{bmatrix*}, \qquad
\mathbf{p}_2 = 
\begin{bmatrix*}[r]
2 \\ 1 \\ 1
\end{bmatrix*}, \qquad
\mathbf{p}_3 = 
\begin{bmatrix*}[r]
1 \\ 0 \\ 0
\end{bmatrix*}.
\end{equation*}
An invertible matrix $P$ can then be formed by stacking these eigenvectors as its columns:
\begin{equation*}
P = \begin{bmatrix*}[r]
 0 & 2 & 1 \\
 1 & 1 & 0 \\
-1 & 1 & 0
\end{bmatrix*}.
\end{equation*}
Its inverse $P^{-1}$ is given by
\begin{equation*}
P^{-1} = \begin{bmatrix*}[r]
0 & \frac{1}{2} & -\frac{1}{2} \\
0 & \frac{1}{2} &  \frac{1}{2} \\
1 & -1 & -1
\end{bmatrix*}.
\end{equation*}
Finally, we can write a diagonalization of $A$ by expressing it as $P D P^{-1}$ and verify that the products are reduced to the original matrix $A$:
\begin{equation*}
\underbrace{\begin{bmatrix*}[r]
1 & 2 & 2 \\
0 &	2 & 1 \\
0 & 1 & 2
\end{bmatrix*}}_{A} = 
\underbrace{\begin{bmatrix*}[r]
 0 & 2 & 1 \\
 1 & 1 & 0 \\
-1 & 1 & 0
\end{bmatrix*}}_{P}
\underbrace{\begin{bmatrix*}[r]
1 & 0 & 0 \\
0 &	3 & 0 \\
0 & 0 & 1
\end{bmatrix*}}_{D}
\underbrace{\begin{bmatrix*}[r]
0 & \frac{1}{2} & -\frac{1}{2} \\
0 & \frac{1}{2} &  \frac{1}{2} \\
1 & -1 & -1
\end{bmatrix*}}_{P^{-1}}.
\end{equation*}

\subsection{Orthogonal diagonalization}		\label{coordo}

The following example illustrates an orthogonal diagonalization for a symmetric matrix, whereby one of the eigenvalues has both algebraic and geometric multiplicities of two, and thus the resulting eigenvectors are not orthogonal. Consider the following $3 \times 3$ matrix $B$:
\begin{equation*}
B = \begin{bmatrix*}[r]
 3 & -4 & -4 \\
-4 &  3 & -4 \\
-4 & -4 &  3
\end{bmatrix*}.
\end{equation*}
This matrix has a characteristic polynomial of $\lambda^3 - 9\lambda^2 - 21\lambda + 245$, which yields three real-valued eigenvalues of $\lambda_1 = -5$ and $\lambda_2 = 7 = \lambda_3$ upon solving the characteristic equation. A diagonal matrix $D$ is given by
\begin{equation*}
D = \begin{bmatrix*}[r]
-5 & 0 & 0 \\
 0 & 7 & 0 \\
 0 & 0 & 7
\end{bmatrix*}.
\end{equation*}
Let $\mathbf{q}_1$, $\mathbf{q}_2$, and $\mathbf{q}_3$ be the eigenvectors corresponding to these eigenvalues, we have
\begin{equation*}
\mathbf{q}_1 = 
\begin{bmatrix*}[r]
1 \\ 1 \\ 1
\end{bmatrix*}, \qquad
\mathbf{q}_2 = 
\begin{bmatrix*}[r]
1 \\ 0 \\ -1
\end{bmatrix*}, \qquad
\mathbf{q}_3 = 
\begin{bmatrix*}[r]
0 \\ 1 \\ -1
\end{bmatrix*}.
\end{equation*}
A common mistake we often encounter is that some students construct an invertible matrix $P$ directly from these eigenvectors without normalizing them, and then find its inverse by simply transposing it, that is, $P^T = P^{-1}$. However, the correct step would be to construct an orthonormal set of eigenvectors that form the column of $P$. We observe that $\mathbf{q}_1$ is orthogonal to both $\mathbf{q}_2$ and $\mathbf{q}_3$. However, $\mathbf{q}_2$ and $\mathbf{q}_3$ are not orthogonal because they originate from the same eigenspace. Only the eigenvectors associated with the distinct eigenvalues are orthogonal. Applying the Gram-Schmidt process to the set $\left\{\mathbf{q}_2, \mathbf{q}_3\right\}$, we can obtain a new set of orthogonal eigenvectors $\left\{\mathbf{q}_2, \mathbf{q}^{\prime}_3\right\}$, where
\begin{equation*}
\mathbf{q}_3^{\prime} = \mathbf{q}_3 - \frac{\left\langle \mathbf{q}_2, \mathbf{q}_3 \right\rangle}{\left\langle \mathbf{q}_2, \mathbf{q}_2 \right\rangle} \mathbf{q}_2 
= \frac{1}{2}\begin{bmatrix*}[r]
-1 \\ 2 \\ -1
\end{bmatrix*}.
\end{equation*}
The angle brackets denote the usual inner (dot) product. By normalizing these vectors, we obtain an orthonormal set of eigenvectors and construct matrix $P$ accordingly. Because $P$ is now an orthogonal matrix, its inverse is simply its transpose, given as follows:
\begin{equation*}
P = \begin{bmatrix*}[r]
\frac{1}{\sqrt{3}} &  \frac{1}{\sqrt{2}} & \frac{1}{\sqrt{6}} \\
\frac{1}{\sqrt{3}} &  0 & -\frac{2}{\sqrt{6}} \\
\frac{1}{\sqrt{3}} & -\frac{1}{\sqrt{2}} & \frac{1}{\sqrt{6}}
\end{bmatrix*}, \qquad 
P^{-1} = P^T = \begin{bmatrix*}[r]
\frac{1}{\sqrt{3}} & \frac{1}{\sqrt{3}} & \frac{1}{\sqrt{3}} \\
\frac{1}{\sqrt{2}} & 0 & -\frac{1}{\sqrt{2}} \\
\frac{1}{\sqrt{6}} & -\frac{2}{\sqrt{6}} & \frac{1}{\sqrt{6}}
\end{bmatrix*}.
\end{equation*}
We can now express an orthogonal diagonalization of $B$ and again confirm that the right-hand side will be reduced to the original matrix upon multiplication and simplification:
\begin{equation*}
\underbrace{\begin{bmatrix*}[r]
 3 & -4 & -4 \\
-4 &  3 & -4 \\
-4 & -4 &  3
\end{bmatrix*}}_{B} = 
\underbrace{\begin{bmatrix*}[r]
\frac{1}{\sqrt{3}} &  \frac{1}{\sqrt{2}} & \frac{1}{\sqrt{6}} \\
\frac{1}{\sqrt{3}} &  0 & -\frac{2}{\sqrt{6}} \\
\frac{1}{\sqrt{3}} & -\frac{1}{\sqrt{2}} & \frac{1}{\sqrt{6}}
\end{bmatrix*}}_{P}
\underbrace{\begin{bmatrix*}[r]
-5 & 0 & 0 \\
 0 & 7 & 0 \\
 0 & 0 & 7
\end{bmatrix*}}_{D}
\underbrace{\begin{bmatrix*}[r]
\frac{1}{\sqrt{3}} & \frac{1}{\sqrt{3}} & \frac{1}{\sqrt{3}} \\
\frac{1}{\sqrt{2}} & 0 & -\frac{1}{\sqrt{2}} \\
\frac{1}{\sqrt{6}} & -\frac{2}{\sqrt{6}} & \frac{1}{\sqrt{6}}
\end{bmatrix*}}_{P^{-T}}.
\end{equation*}

\subsection{SVD}		\label{svdsubsec23}

Whereas matrix diagonalization and orthogonal diagonalization can be implemented in a square matrix, SVD encompasses both square and rectangular matrices. Let $M$ be an $m \times n$ real-valued matrix with rank $r$, then an SVD of $M$ is given by $M = U \Sigma V^T$, where $U$ is an $m \times m$ orthogonal matrix, $\Sigma$ is an $m \times n$ rectangular ``diagonal'' matrix with non-negative real numbers on the diagonal, and $V$ is an $n \times n$  is orthogonal matrix. The columns of $U$ and $V$ in such a decomposition are called the left and right singular vectors of $M$, respectively. For complex-valued matrices, $U$ and $V$ are complex unitary matrices, and instead of $V^T$, we would have the conjugate transpose of $V$, that is, $V^{\ast}$. In general, $\Sigma$ is composed of matrix blocks and admits the following form:
\begin{equation*}
\Sigma = \begin{bmatrix*}[c]
D & 0 \\ 0 & 0
\end{bmatrix*},
\end{equation*}
where $D$ is an $r \times r$ diagonal matrix, for $r \leq \min\{m, n\}$. The number of (zero) rows and columns in the second row and column of the block are $(m - r)$ and $(n - r)$, respectively. For $r = m$, $r = n$, or $r = m = n$, some or all zero blocks would disappear. Furthermore, the diagonal entries of $D$ are the first $r$ singular values of $M$, that is, $\sigma_1 \geq \sigma_2 \geq \cdots \geq \sigma_r > 0$.

Similar to obtaining matrix diagonalization and orthogonal diagonalization, constructing an SVD of a matrix requires a step-by-step algorithm and an understanding of (orthogonal) diagonalization. For square matrices with full rank, many learners usually encounter no difficulty when constructing an SVD of a matrix, although one must be careful when calculating eigenvalues and the corresponding eigenvectors of a symmetric matrix $M^T M$ instead of the original matrix $M$. However, for rectangular matrices or those with rank deficiencies, many students often could not complete the construction of an SVD, partly because one or more singular values can be zero. The following example illustrates the construction of an SVD of a square matrix with rank deficiency. Let
\begin{equation*}
M = \begin{bmatrix*}[r]
3 & -1 \\ -6 & 2
\end{bmatrix*}
\end{equation*}
be a matrix with rank $r = 1$, then the eigenvalues of $A = M^T M$ are given by $\lambda_1 = 50$ and $\lambda_2 = 0$, which gives only one singular value $\sigma_1 = 5 \sqrt{2}$ because the second eigenvalue is zero. Thus, $D = \sigma_1$ and matrix $\Sigma$ contains one block of zero in each row and column below and to the right of $D$, respectively, given as follows: 
\begin{equation*}
\Sigma = \begin{bmatrix*}[c]
5 \sqrt{2}	& 0 \\ 0 & 0
\end{bmatrix*}.
\end{equation*}
By finding the corresponding eigenvectors of $\lambda_1$ and $\lambda_2$, we obtain $\mathbf{v}_1$ and $\mathbf{v}_2$ and set up an orthogonal matrix $V$:
\begin{equation*}
V = \begin{bmatrix*}[r]
3/\sqrt{10}	& 1/\sqrt{10} \\ -1/\sqrt{10} & 3/\sqrt{10}
\end{bmatrix*}.
\end{equation*}
The first column of matrix $U$, i.e., $\mathbf{u}_1$, can be calculated using a formula
\begin{equation*}
\mathbf{u}_1 = \frac{1}{\sigma_1} M \mathbf{v}_1 = \frac{1}{5 \sqrt{2}} 
\begin{bmatrix*}[r]
3 & -1 \\ -6 & 2
\end{bmatrix*}
\frac{1}{\sqrt{10}}
\begin{bmatrix*}[r]
3 \\ -1
\end{bmatrix*} = 
\frac{1}{\sqrt{5}}
\begin{bmatrix*}[r]
1 \\ -2
\end{bmatrix*}.
\end{equation*}
The second column of matrix $U$, i.e., $\mathbf{u}_2$, must be orthonormal to $\mathbf{u}_1$:
\begin{equation*}
\mathbf{u}_2 = \frac{1}{\sqrt{5}}
\begin{bmatrix*}[r]
2 \\ 1
\end{bmatrix*}.
\end{equation*}
An SVD of $M$ is given as follows:
\begin{equation*}
\underbrace{\begin{bmatrix*}[r]
3 & -1 \\ -6 & 2
\end{bmatrix*}}_{M} = 
\underbrace{\begin{bmatrix*}[r]
2/\sqrt{5} & 1/\sqrt{5} \\ 1/\sqrt{5} & -2/\sqrt{5}
\end{bmatrix*}}_{U}
\underbrace{\begin{bmatrix*}[c]
5 \sqrt{2}	& 0 \\ 0 & 0
\end{bmatrix*}}_{\Sigma}
\underbrace{\begin{bmatrix*}[r]
3/\sqrt{10}	& -1/\sqrt{10} \\ 1/\sqrt{10} & 3/\sqrt{10}
\end{bmatrix*}}_{V^T}.
\end{equation*}

\section{(Orthogonal) diagonalization and SVD using SageMath}		\label{sagemathsection}

To use SageMath online, we can utilize Sage Cell Server or SageMathCell, accessible online at the URL \url{https://sagecell.sagemath.org/}. In addition to an open-source, scalable, and easy-to-use web interface for SageMath, this cell server also allows embedding Sage computations into any webpage. Without loading up any program, this is one way to conduct one-off computations using SageMath, with the idea of accessing the computation in the cloud as simply as possible, as long as one has an internet connection. 

\subsection{Matrix diagonalization using SageMath}

The following provides the SageMath command for matrix diagonalization discussed in Subsection~\ref{madi21}.
\begin{Verbatim}
A=matrix([[1,2,2],[0,2,1],[0,1,2]])
print("Matrix A = ")
print(A)
print()
print("Characteristic polynomial of A: p(x) =", A.charpoly())
print()
print("Eigenvalues of A =", A.eigenvalues())
print()
print("Eigenvalue, eigenvector, and geometric multiplicity:", )
print(A.eigenvectors_right())
print()
lambda1 = A.eigenvalues()[0]
lambda2 = A.eigenvalues()[1]
lambda3 = A.eigenvalues()[2]
print("Extracting eigenvalues:")
print("lambda1 =", lambda1)
print("lambda2 =", lambda2)
print("lambda3 =", lambda3)
print()
D=diagonal_matrix(A.eigenvalues())
print("Diagonal matrix D = ")
print(D)
print()
r1=A.eigenvectors_right()[0][1]
p1=2*r1[0]
r2=A.eigenvectors_right()[1][1]
p2=r2[0]
r3=A.eigenvectors_right()[1][1]
p3=r3[1]
print("Extracting eigenvectors:")
print("p1 =", p1)
print("p2 =", p2)
print("p3 =", p3)
print()
P=matrix([p1,p2,p3]).transpose()
print("Invertible matrix P = ")
print(P)
print()
print("Inverse of P, P^(-1) = ")
print(P.inverse())
print()
print("Calculate PDP^(-1) = ")
print(P*D*P.inverse())
print("= A")
\end{Verbatim}
Once we know the eigenvalues and their corresponding eigenvectors, the construction of diagonal matrix $D$ and invertible matrix $P$ follows immediately. Indeed, this process must be performed manually, but this is not a difficult step. Furthermore, finding the inverse of $P$ follows naturally using the Sage's command ``\verb*|P.inverse()|''. The corresponding output of the previous lines of code is given as follows:
\begin{Verbatim}
Matrix A = 
[1 2 2]
[0 2 1]
[0 1 2]

Characteristic polynomial of A: p(x) = x^3 - 5*x^2 + 7*x - 3

Eigenvalues of A = [3, 1, 1]

Eigenvalue, eigenvector, and geometric multiplicity:
[(3, [(1, 1/2, 1/2)], 1), (1, [(1, 0, 0),(0, 1, -1)], 2)]

Extracting eigenvalues:
lambda1 = 3
lambda2 = 1
lambda3 = 1

Diagonal matrix D = 
[3 0 0]
[0 1 0]
[0 0 1]

Extracting eigenvectors:
p1 = (2, 1, 1)
p2 = (1, 0, 0)
p3 = (0, 1, -1)

Invertible matrix P = 
[ 2  1  0]
[ 1  0  1]
[ 1  0 -1]

Inverse of P, P^(-1) = 
[   0  1/2  1/2]
[   1   -1   -1]
[   0  1/2 -1/2]

Calculate PDP^(-1) = 
[1 2 2]
[0 2 1]
[0 1 2]
= A
\end{Verbatim}
In the final step, we confirm that the diagonalization is correct as the product returns to the original matrix $A$.

\subsection{Orthogonal diagonalization using SageMath}

The following commands illustrate the SageMath codes for an example of an orthogonal diagonalization of the matrix considered in Subsection~\ref{coordo}.
\begin{Verbatim}
A=matrix([[3,-4,-4],[-4,3,-4],[-4,-4,3]])
print("Matrix A = ")
print(A)
print()
print("Characteristic polynomial of A: p(x) =", A.charpoly())
print()
print("Eigenvalues of A =", A.eigenvalues())
print()
lambda1 = A.eigenvalues()[0]
lambda2 = A.eigenvalues()[1]
lambda3 = A.eigenvalues()[2]
print("Extracting eigenvalues:")
print("lambda1 =", lambda1)
print("lambda2 =", lambda2)
print("lambda3 =", lambda3)
print()
D=diagonal_matrix(A.eigenvalues())
print("Diagonal matrix D = ")
print(D)
print()
print("Eigenvalue, eigenvector, and geometric multiplicity:", )
print(A.eigenvectors_right())
print()
q1=A.eigenvectors_right()[0][1][0]
q2=A.eigenvectors_right()[1][1][0]
q3=A.eigenvectors_right()[1][1][1]
print("Extracting eigenvectors:")
print("q1 =", q1)
print("q2 =", q2)
print("q3 =", q3)
print()
print("Checking orthogonality:")
print("<q1,q2> =", q1.inner_product(q2))
print("<q1,q3> =", q1.inner_product(q3))
print("<q2,q3> =", q2.inner_product(q3))
print()
print("q2 and q3 are not orthogonal; apply the Gram-Schmidt process:")
print()
q31 = q3 - q2.inner_product(q3)/q2.inner_product(q2)*q2
print("q3' =", q31)
print()
print("Check that now q2 and q3' are orthogonal: <q2,q3'> =", q2.inner_product(q31))
print()
print("Normalize all the orthogonal eigenvectors:")
print()
p1 = q1/sqrt(q1.inner_product(q1))
print("p1 =", p1)
print()
p2 = q2/sqrt(q2.inner_product(q2))
print("p2 =", p2)
print()
q32 = 2*q31
p3 = -q32/sqrt(q32.inner_product(q32))
print("p3 =", p3)
print()
print("Construct an invertible and orthogonal matrix P:")
P = matrix([p1,p2,p3]).transpose()
print("P = ")
print(P)
print()
print("Inverse of P, P^(-1) = ")
print(P.inverse())
print()
print("Calculate PDP^(-1) = ")
print(P*D*P.inverse())
print("= A")
print()
print("Transpose of P, P^T = ")
print(P.transpose())
print()
print("Calculate PDP^T = ")
print(P*D*P.transpose())
print("= A")
\end{Verbatim}
The output that we obtained is given as follows:
\begin{Verbatim}
Matrix A = 
[ 3 -4 -4]
[-4  3 -4]
[-4 -4  3]

Characteristic polynomial of A: p(x) = x^3 - 9*x^2 - 21*x + 245

Eigenvalues of A = [-5, 7, 7]

Extracting eigenvalues:
lambda1 = -5
lambda2 = 7
lambda3 = 7

Diagonal matrix D = 
[-5  0  0]
[ 0  7  0]
[ 0  0  7]

Eigenvalue, eigenvector, and geometric multiplicity:
[(-5, [(1, 1, 1)], 1), (7, [(1, 0, -1),(0, 1, -1)], 2)]

Extracting eigenvectors:
q1 = (1, 1, 1)
q2 = (1, 0, -1)
q3 = (0, 1, -1)

Checking orthogonality:
<q1,q2> = 0
<q1,q3> = 0
<q2,q3> = 1

q2 and q3 are not orthogonal; apply the Gram-Schmidt process:

q3' = (-1/2, 1, -1/2)

Check that now q2 and q3' are orthogonal: <q2,q3'> = 0

Normalize all the orthogonal eigenvectors:

p1 = (1/3*sqrt(3), 1/3*sqrt(3), 1/3*sqrt(3))

p2 = (1/2*sqrt(2), 0, -1/2*sqrt(2))

p3 = (1/6*sqrt(6), -1/3*sqrt(6), 1/6*sqrt(6))

Construct an invertible and orthogonal matrix P:
P = 
[ 1/3*sqrt(3)  1/2*sqrt(2)  1/6*sqrt(6)]
[ 1/3*sqrt(3)            0 -1/3*sqrt(6)]
[ 1/3*sqrt(3) -1/2*sqrt(2)  1/6*sqrt(6)]

Inverse of P, P^(-1) = 
[ 1/3*sqrt(3)  1/3*sqrt(3)  1/3*sqrt(3)]
[ 1/2*sqrt(2)            0 -1/2*sqrt(2)]
[ 1/6*sqrt(6) -1/3*sqrt(6)  1/6*sqrt(6)]

Calculate PDP^(-1) = 
[ 3 -4 -4]
[-4  3 -4]
[-4 -4  3]
= A

Transpose of P, P^T = 
[ 1/3*sqrt(3)  1/3*sqrt(3)  1/3*sqrt(3)]
[ 1/2*sqrt(2)            0 -1/2*sqrt(2)]
[ 1/6*sqrt(6) -1/3*sqrt(6)  1/6*sqrt(6)]

Calculate PDP^T = 
[ 3 -4 -4]
[-4  3 -4]
[-4 -4  3]
= A
\end{Verbatim}
Because $P^{-1} = P^T$, it does not come as a surprise that the calculation of both $PDP^{-1}$ and $PDP^T$ is reduced to the original matrix $A$. However, in this example, although we can directly construct a diagonal matrix $D$, we should apply the Gram-Schmidt process to construct an orthogonal and invertible matrix $P$.

\subsection{SVD using SageMath}

The following provides the SageMath command for constructing an SVD of the matrix discussed in Subsection~\ref{svdsubsec23}.
\begin{Verbatim}
M = matrix([ [2,-2], [-3,-4], [-4,-3] ])
print("Matrix M = ")
print(M)
print()
A = M.transpose()*M
print("Matrix M^T M = ")
print(A)
print()
print("Eigenvalues of M^T M =", A.eigenvalues())
print()
print("Eigenvalue, eigenvector, and geometric multiplicity:", )
print(A.eigenvectors_right())
print()
sigma1 = sqrt(A.eigenvalues()[0])
sigma2 = sqrt(A.eigenvalues()[1])
print("Singular value of M:")
print("sigma1 =", sigma1)
print("sigma2 =", sigma2)
print()
Sigma = matrix([[sigma1, 0], [0, sigma2], [0,0]])
print("Matrix Sigma =")
print(Sigma)
print()
w1 = A.eigenvectors_right()[0][1][0]
v1 = w1/w1.norm()
w2 = A.eigenvectors_right()[1][1][0]
v2 = w2/w2.norm()
print("Eigenvectors of M^T M:")
print("v1 =", v1)
print("v2 =", v2)
print()
print("Matrix V =")
V = matrix([ v1, v2 ])
print(V)
print()
u1 = 1/sigma1*M*v1
u2 = 1/sigma2*M*v2
u3 = u1.cross_product(u2)
print("Left singular vectors of M:")
print("u1 =", u1)
print("u2 =", u2)
print("u3 =", u3)
print()
print("Checking orthogonality")
print("u1.u2 =", u1.dot_product(u2))
print("u1.u3 =", u1.dot_product(u3))
print("u2.u3 =", u3.dot_product(u2))
print()
UT = matrix([u1, u2, u3])
U = UT.transpose()
print("Matrix U =")
print(U)
print()
print("SVD of M = U*Sigma*V^T =")
print(U*Sigma*V.transpose())
print("= M")
\end{Verbatim}
The output is given as follows.
\begin{Verbatim}
Matrix M = 
[ 2 -2]
[-3 -4]
[-4 -3]

Matrix M^T M = 
[29 20]
[20 29]

Eigenvalues of M^T M = [49, 9]

Eigenvalue, eigenvector, and geometric multiplicity:
[(49, [(1, 1)], 1), (9, [(1, -1)], 1)]

Singular value of M:
sigma1 = 7
sigma2 = 3

Matrix Sigma =
[7 0]
[0 3]
[0 0]

Eigenvectors of M^T M:
v1 = (1/2*sqrt(2), 1/2*sqrt(2))
v2 = (1/2*sqrt(2), -1/2*sqrt(2))

Matrix V =
[ 1/2*sqrt(2)  1/2*sqrt(2)]
[ 1/2*sqrt(2) -1/2*sqrt(2)]

Left singular vectors of M:
u1 = (0, -1/2*sqrt(2), -1/2*sqrt(2))
u2 = (2/3*sqrt(2), 1/6*sqrt(2), -1/6*sqrt(2))
u3 = (1/3, -2/3, 2/3)

Checking orthogonality
u1.u2 = 0
u1.u3 = 0
u2.u3 = 0

Matrix U =
[           0  2/3*sqrt(2)          1/3]
[-1/2*sqrt(2)  1/6*sqrt(2)         -2/3]
[-1/2*sqrt(2) -1/6*sqrt(2)          2/3]

SVD of M = U*Sigma*V^T =
[ 2 -2]
[-3 -4]
[-4 -3]
= M
\end{Verbatim}
We confirm that the construction of our SVD is correct because the product of each matrix returns to the original matrix $M$, that is, $M = U \Sigma V^T$. In this example, both vectors $\mathbf{u}_1$ and $\mathbf{u}_2$ are calculated using a formula that involves a singular value, that is, 
\begin{equation*}
\mathbf{u}_i = \frac{1}{\sigma_i} M \mathbf{v}_i, \qquad \qquad \text{for} \quad i = 1, 2.
\end{equation*}
However, because the third singular value $\sigma_3$ does not exist, we cannot use this formula but calculate $\mathbf{u}_3$ either using a dot or cross product to ensure its orthogonality with $\mathbf{u}_1$ and $\mathbf{u}_2$. We applied the latter to this example. 

We observed that some students did not know how to find the third vector $\mathbf{u}_3$. Very often, they simply ignored and abandoned it entirely, and construct a factorization using the first two vectors in matrix $U$ instead, that is, $\mathbf{u}_1$ and $\mathbf{u}_2$. In this situation, they confused between ``reduced SVD'' and ``(complete, non-reduced) SVD''. The former depends on the rank of $M$, i.e., $r$, and instead of employing $\Sigma$, we admit diagonal matrix $D$ instead. Thus, the reduced SVD for $M$ in the example considered in this subsection is given by
\begin{equation*}
\underbrace{\begin{bmatrix*}[r]
2 & -2 \\ -3 & -4 \\ -4 & -3
\end{bmatrix*}}_{M} = 
\underbrace{\begin{bmatrix*}[r]
 0 & 2\sqrt{2}/3 \\ 
-1/\sqrt{2} &  \sqrt{2}/6 \\
-1/\sqrt{2} & -\sqrt{2}/6
\end{bmatrix*}}_{U_2}
\underbrace{\begin{bmatrix*}[c]
7 & 0 \\ 0 & 3
\end{bmatrix*}}_{D}
\underbrace{\begin{bmatrix*}[r]
1/\sqrt{2}	& 1/\sqrt{2} \\ 1/\sqrt{2} & -1/\sqrt{2}
\end{bmatrix*}}_{V_2^T}.
\end{equation*}
This reduced SVD is beneficial for calculating the pseudoinverse $M^{+}$ (also called the Moore-Penrose inverse) of the matrix $M$. For the particular $M$ considered in this subsection, its pseudoinverse is given by
\begin{equation*}
\underbrace{\frac{1}{63}\begin{bmatrix*}[r]
14 & -1 & -8 \\ -14 & -8 & -1
\end{bmatrix*}}_{M^+} = 
\underbrace{\begin{bmatrix*}[r]
1/\sqrt{2}	& 1/\sqrt{2} \\ 1/\sqrt{2} & -1/\sqrt{2}
\end{bmatrix*}}_{V_2}
\underbrace{\begin{bmatrix*}[c]
1/7 & 0 \\ 0 & 1/3
\end{bmatrix*}}_{D^{-1}}
\underbrace{\begin{bmatrix*}[r]
0 &	-1/\sqrt{2} & -1/\sqrt{2} \\
2\sqrt{2}/3 &  \sqrt{2}/6 & -\sqrt{2}/6
\end{bmatrix*}}_{U_2^T}.
\end{equation*}
SageMath can certainly compute the pseudoinverse of a matrix in a matter of seconds using the command \verb*|print(M.pseudoinverse())|. However, it is important to exercise discretion and only use this step for checking the final calculation rather than as a substitute for an appropriate algorithm to find it.

\section{(Orthogonal) diagonalization and SVD in ChatGPT}		\label{chatgptsection}

We used the version from March 14, 2023, for any responses generated by ChatGPT in this study.

\subsection{Matrix diagonalization in ChatGPT}

When asking the chatbot to provide an example of a matrix diagonalization without specifying any size of the matrix, ChatGPT provided an answer by selecting a simple $2 \times 2$ matrix, as illustrated in the following example. We asked ``Provide an example of a matrix diagonalization.'' ChatGPT responded:\\
{\ttfamily \noindent
Sure, here's an example of matrix diagonalization: \\
Suppose we have the following matrix $A$:
\begin{equation*}
\begin{bmatrix*}[c]
3 & 1 \\ 0 & 2
\end{bmatrix*}.
\end{equation*}
We want to diagonalize $A$, which means finding a diagonal matrix $D$ and an invertible\\ matrix $P$ such that $A = PDP^{-1}$, where $D$ has the eigenvalues of $A$ along the diagonal. \\
To diagonalize $A$, we first find the eigenvalues of $A$. The characteristic polynomial\\ of $A$ is:
\begin{align*}
\det\left(A - \lambda I \right) 
&= \det \left(
\begin{matrix*}[c]
3 - \lambda & 1 \\
0  & 2 - \lambda
\end{matrix*}
\right) \\
&= \left(3 - \lambda \right) \left(2 - \lambda \right) = 6 - 5 \lambda + \lambda^2.	
\end{align*}
Solving for the roots of the characteristic polynomial, we get $\lambda = 3$ and $\lambda = 2$. These are the eigenvalues of $A$. \\
Next, we find the eigenvectors corresponding to each eigenvalue. For $\lambda = 3$, we solve the equation $\left(A - 3I \right) \mathbf{x} = \mathbf{0}$:
\begin{equation*}
\begin{bmatrix*}[r]
0 & 1 \\ 0 & -1
\end{bmatrix*}
\begin{bmatrix*}[r]
x_1 \\ x_2
\end{bmatrix*} = 
\begin{bmatrix*}[r]
0 \\ 0 
\end{bmatrix*}.
\end{equation*}
This gives us the equation $x_2 = 0$ and $x_1$ can be any nonzero value. So an eigenvector corresponding to $\lambda = 3$ is $\left[1 \quad 0 \right]^T$. \\
For $\lambda = 2$, we solve the equation $\left(A - 2I \right) \mathbf{x} = \mathbf{0}$:
\begin{equation*}
\begin{bmatrix*}[r]
1 & 1 \\ 0 & 0
\end{bmatrix*}
\begin{bmatrix*}[r]
x_1 \\ x_2
\end{bmatrix*} = 
\begin{bmatrix*}[r]
0 \\ 0 
\end{bmatrix*}.
\end{equation*}
This gives us the equation $x_1 + x_2 = 0$, or $x_1 = -x_2$. So an eigenvector corresponding to $\lambda = 2$ is $\left[-1 \quad 1 \right]^T$. \\
We can normalize these eigenvectors to get:
\begin{align*}
\mathbf{v}_1 &= \frac{\left[1 \quad 0 \right]^T}{\sqrt{1^2 + 0^2}} = \left[1 \quad 0 \right]^T, \\
\mathbf{v}_2 &= \frac{\left[-1 \quad 1 \right]^T}{\sqrt{1^2 + 1^2}} = \left[-1/\sqrt{2} \quad 1/\sqrt{2} \right]^T.
\end{align*}
We can then form the matrix $P$ by stacking the eigenvectors as columns:
\begin{equation*}
P = \begin{bmatrix*}[r]
1 & -1/\sqrt{2} \\
0 &  1/\sqrt{2}
\end{bmatrix*}.
\end{equation*}
We can also form the diagonal matrix $D$ by putting the eigenvalues along the diagonal:
\begin{equation*}
D = \begin{bmatrix*}[r]
3 & 0 \\
0 & 2
\end{bmatrix*}.
\end{equation*}
Finally, we can verify that $A = PDP^{-1}$:
\begin{equation*}
P D P^{-1} = 
\begin{bmatrix*}[r]
1 & -1/\sqrt{2} \\
0 &  1/\sqrt{2}
\end{bmatrix*}
\begin{bmatrix*}[r]
3 & 0 \\
0 & 2
\end{bmatrix*}
\begin{bmatrix*}[r]
 1/\sqrt{2} & 1/\sqrt{2} \\
-1/\sqrt{2}	& 1/\sqrt{2}
\end{bmatrix*} = 
\begin{bmatrix*}[r]
3 & 1 \\
0 & 2
\end{bmatrix*},
\end{equation*}
which is equal to $A$. So we have successfully diagonalized $A$.}

As we observe in this example, ChatGPT came up with an excellent diagonalization example with a step-by-step explanation that can help readers better understand the proper algorithm. Although the calculation for finding $p^{-1}$ was not shown, the machine assumed that the readers have a working knowledge of it, particularly for an easy, simple $2 \time 2$ matrix where the computation is relatively straightforward. As some examples in the following subsection illustrate, ChatGPT encountered some difficulty in (orthogonally) diagonalizing a larger matrix size, such as $3 \times 3$.

\subsection{Orthogonal diagonalization in ChatGPT}

When we asked ChatGPT to produce an example of diagonalization of a $3 \times 3$ symmetric matrix, it gave the following response:\\
{\ttfamily
Sure, let's consider the symmetric matrix:
\begin{equation*}
A = \begin{bmatrix*}[r]
1 & 2 & 3 \\
2 & 5 & 6 \\
3 & 6 & 9
\end{bmatrix*}.
\end{equation*}
To diagonalize this matrix, we need to find a matrix $P$ and a diagonal matrix $D$ such that:
\begin{equation*}
A = P D P^T,
\end{equation*}
where $P^T$ denotes the transpose of $P$. First, we need to find the eigenvalues of $A$ by solving the characteristic equation:
\begin{equation*}
\det(A - \lambda I) = 0,
\end{equation*}
where $\lambda$ is the eigenvalue and $I$ is the identity matrix. This gives us:
\begin{align*}
&\det \begin{bmatrix*}[c]
1 - \lambda & 2 		  & 3 \\
		2 	& 5 - \lambda & 6 \\
		3 	& 6 		  & 9 - \lambda
\end{bmatrix*} = 0\\
&= \left(1 - \lambda\right) \left[\left(5 - \lambda \right) \left(9 - \lambda \right) - 6\ast6 \right] -2 \left[2*\left(9 - \lambda \right) - 6* 3 \right] + 3 \left[2*6 - {\color{red} 5}*3 \right] \\
&= \left(1- \lambda \right) \left(\lambda^2 - {\color{red} 15} \lambda + {\color{red} 36} \right) + 2 {\color{red} \left(3 \lambda - 12 \right)} - 3 {\color{red} \left(6 - 15 \right)} \\
&= \lambda^3 - 15 \lambda^2 + {\color{red} 54 \lambda - 36}
\end{align*}
}

We observe that ChatGPT made calculation mistakes starting with the second expression, where it wrote $5$ instead of $5 - \lambda$. The subsequent terms, including the final expression for the characteristic equation, were undoubtedly incorrect. They are represented by the red terms. The correct characteristic equation was $\lambda^3 - 15 \lambda^2 + 10 \lambda = 0$. Let us continue with what ChatGPT offers after this computational step.\\

{\ttfamily \noindent
This polynomial has a repeated root at $\lambda = 0$, and distinct roots at $\lambda = 3$ and $\lambda = 12$. Therefore, the eigenvalues of $A$ are $\lambda_1 = 0$, $\lambda_2 = 3$, and $\lambda_3 = 12$. Next, we need to find the eigenvectors of $A$ for each eigenvalue. For $\lambda_1 = 0$, we have:
\begin{equation*}
A - \lambda_1 I = A = 
\begin{bmatrix*}[c]
1 & 2 & 3 \\
2 & 5 & 6 \\
3 & 6 & 9 \\
\end{bmatrix*}.
\end{equation*}
Reducing this matrix to row echelon form, we get:
\begin{equation*}
\begin{bmatrix*}[c]
1 & 2 & 3 \\
0 & 1 & {\color{red} 1} \\
0 & 0 & 0 \\
\end{bmatrix*},
\end{equation*}
which gives us the eigenvector $v_1 = [-2, \ 1, \ 0]$. For $\lambda_2 = 3$, we have:
\begin{equation*}
A - \lambda_2 I =  
\begin{bmatrix*}[r]
-2 & 2 & 3 \\
 2 & 2 & 6 \\
 3 & 6 & 6 \\
\end{bmatrix*}.
\end{equation*}
Reducing this matrix to row echelon form, we get:
\begin{equation*}
\begin{bmatrix*}[c]
1 & 0 & 1 \\
0 & 1 & 2 \\
0 & 0 & 0 \\
\end{bmatrix*},
\end{equation*}
which gives us the eigenvector $v_2 = [-1, \ -2, \ 1]$. For $\lambda_3 = 12$, we have:
\begin{equation*}
A - \lambda_3 I =  
\begin{bmatrix*}[r]
-11 &  2 &  3 \\
  2 & -7 &  6 \\
  3 &  6 & -3 \\
\end{bmatrix*}.
\end{equation*}
Reducing this matrix to row echelon form, we get:
\begin{equation*}
\begin{bmatrix*}[r]
1 & 0 & -1 \\
0 & 1 & -2 \\
0 & 0 & 0 \\
\end{bmatrix*},
\end{equation*}
which gives us the eigenvector $v_3 = [1, \ 2, \ {\color{red}1]}$.
}

ChatGPT stopped at this stage, even without completing the final component for the third eigenvector, as indicated by the red color. In other instances, ChatGPT could continue the process by constructing the diagonal matrix $D$, an invertible matrix $P$, and finalizing the answer with ``{\ttfamily Finally, we can check that $A = PDP^{-1}$.}'', this did not occur in this example.

Interestingly, however, ChatGPT obtained one correct eigenvalue, that is, $\lambda_1 = 0$, although one entry in the row echelon form matrix is incorrect, that is, $a_{23} = 0$ instead of 1. The correct associated eigenvector was $\mathbf{v}_1 = \left(3, \ 0, \ -1 \right)$. Furthermore, contrary to other positive integer eigenvalues obtained by ChatGPT, $\lambda_2$ and $\lambda_3$ of matrix $A$ in this example are not integers, but irrational numbers instead, although both are still positive. Arranged ascendingly, they are
\begin{equation*}
\lambda_2 = \frac{1}{2} \left(15 - \sqrt{185} \right), \qquad \qquad \text{and} \qquad \qquad \lambda_3 = \frac{1}{2} \left(15 + \sqrt{185} \right).
\end{equation*}
A possible choice for the corresponding eigenvectors would be
\begin{equation*}
\mathbf{v}_2 = \left(\frac{1}{3}, \ \frac{2 \left(-15 + \sqrt{185}\right) }{3 \left(-13 + \sqrt{185}\right)}, \ 1 \right), \qquad \qquad \text{and} \qquad \qquad
\mathbf{v}_3 = \left(\frac{1}{3}, \ \frac{2 \left( 15 + \sqrt{185}\right) }{3 \left( 13 + \sqrt{185}\right)}, \ 1 \right).
\end{equation*}
If we use the eigenvalues and associated eigenvectors that ChatGPT delivered them out to construct the original matrix, the resulting matrix is neither symmetric nor has all integer entries. Instead, we have the following matrix:
\begin{equation*}
P D P^{-1} = \frac{1}{10}
\begin{bmatrix*}[r]
15 & 30 & 45 \\
30 & 60 & 90 \\
 9 & 18 & 75\\
\end{bmatrix*} = 
\begin{bmatrix*}[c]
\frac{3}{2} & 3 & \frac{9}{2}\\
 3 & 6 & 9 \\
\frac{9}{10} & \frac{9}{5} & \frac{15}{2} \\
\end{bmatrix*}.
\end{equation*}

To further test its capability, we inquired ChatGPT to orthogonally diagonalize a particular $3 \time 3$ symmetric matrix. This is the same matrix $B$ that we explore in Subsection~\ref{coordo}:
\begin{equation*}
B= \begin{bmatrix*}[r]
 3 & -4 & -4 \\
-4 &  3 & -4 \\
-4 & -4 &  3
\end{bmatrix*}.
\end{equation*}
Unfortunately, ChatGPT provided an incorrect answer for this matrix diagonalization. Its computation for the characteristic polynomial was erroneous: $\lambda^3 - 9 \lambda^2 + 24 \lambda$, which gives a false set of eigenvalues: $\lambda_1 = 0$, $\lambda_2 = 3$, and $\lambda_3 = 6$. The calculation was continued to seek the corresponding eigenvectors $\mathbf{v}_1$ and $\mathbf{v}_2$ but was not able to complete the calculation for $\mathbf{v}_3$. The computation stopped abruptly with some missing entries in the matrix equation. By further analyzing each step, we discovered that neither the row-reduction process of the augmented matrix for the first and second eigenvalues nor the obtained eigenvector was correct. Requesting ChatGPT to ``regenerate response'' does not really help either, as it again yielded wrong and incomplete answers. From these examples, we observe that the current version of ChatGPT does not seem ready for the (orthogonal) diagonalization of a simple matrix beyond $2 \times 2$.

\subsection{SVD in ChatGPT}
We asked ChatGPT to ``provide an example of a singular value decomposition.'' Interestingly, ChatGPT provided a relatively advanced example of a $3 \times 3$ matrix and utilized a matrix algebra library \texttt{numpy} in Python, where SageMath also employs it as its primary language solver. The chosen matrix is given by
\begin{equation*}
A = \begin{bmatrix*}[r]
 4 & 11 & 14 \\
-1 &  1 & 17 \\
 0 & 12 &  5
\end{bmatrix*}.
\end{equation*}
In the previous two cases, ChatGPT explained the algorithm to perform an (orthogonal) diagonalization of a matrix. In this case, it directly provided a Python code and its corresponding output. Furthermore, although the original matrix contains integer entries, the matrices in its SVD consist of real numbers, which can be difficult to check manually.
\begin{figure}[H]
\centering
\includegraphics[width=0.49\textwidth]{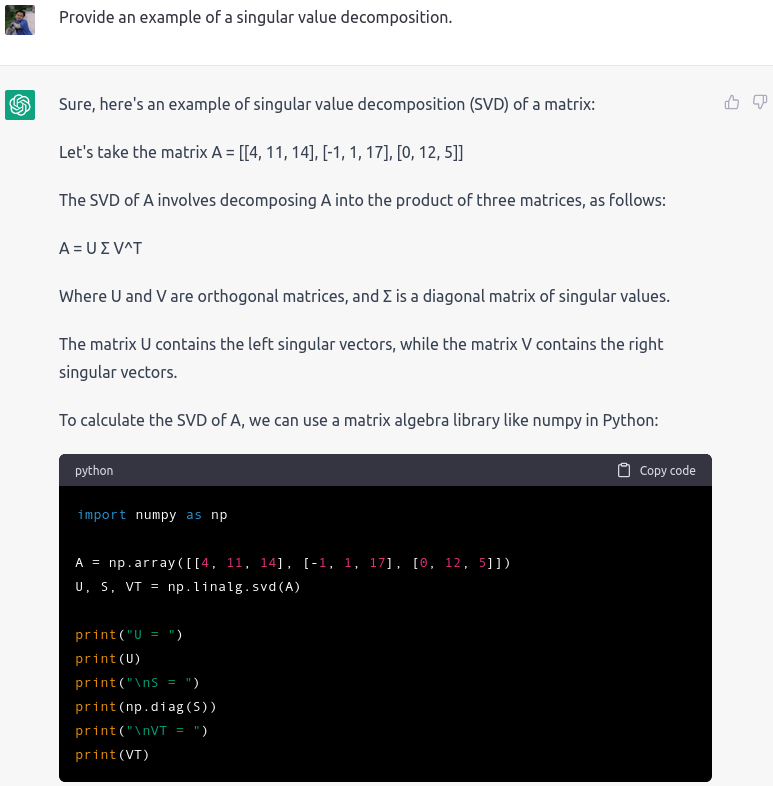} \hspace{0.05cm}
\includegraphics[width=0.49\textwidth]{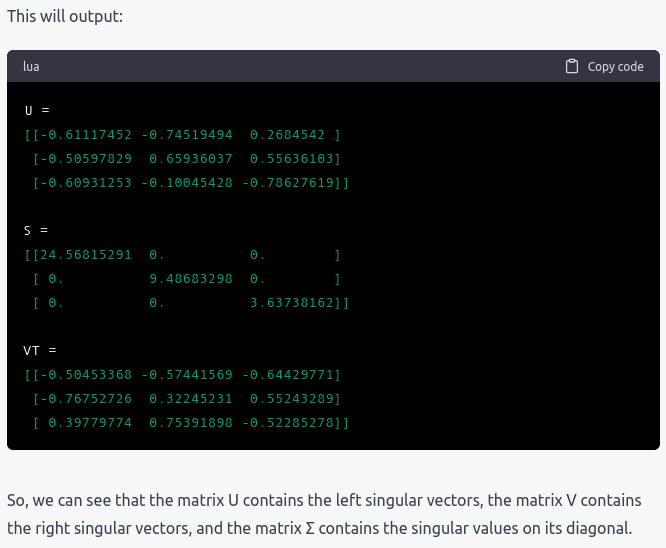}
\caption{An output of ChatGPT to the inquiry ``Provide an example of a singular value decomposition.''} \label{svd1chat}
\end{figure}

Unfortunately, ChatGPT got it wrong for its own given example, as shown in Figure~\ref{svd1chat}. It dispensed singular values of $\sigma_1 = 24.57$, $\sigma_2 = 9.49$, and $\sigma_3 = 3.63$, whereas the correct singular values are $\sigma_1 = 25.58$, $\sigma_2 = 11.36$, and $\sigma_3 = 3.13$. Although the error for the first singular value is less than 4\%, the errors for the second and third singular values are approximately 16\%, which is too large to be considered acceptable. Using the SageMath commands \verb*|A=matrix(RDF,[[4,11,14],[-1,1,17],[0,12,5]])| and \verb*|A.SVD()|, we obtain the following output
\begin{Verbatim}
(
[-0.7020788202506024 0.21713691817311598 -0.6781864706124737]
[-0.5829773893610252 -0.7221727564849129  0.3722954113130941]
[-0.4089287146304667  0.6567481012970376  0.6336081105805151],

[25.575254913906992                0.0                0.0]
[               0.0 11.363233131054077                0.0]
[               0.0                0.0 3.1278217497152005],

[-0.08701136700816689  0.13998836518016763  -0.9863226042353253]
[ -0.5166317611363659   0.8401931429966347  0.16482446980522636]
[ -0.8517749969308255  -0.5239071865088518 0.000783918281394777]
),
\end{Verbatim}
where the first, second, and third matrices correspond to $U$, $\Sigma$, and $V^T$, respectively. Here, RDF stands for Real Double Field, an approximation to a real number using double precision floating point numbers. 

\begin{figure}[h]
\begin{center}
\includegraphics[width=0.49\textwidth]{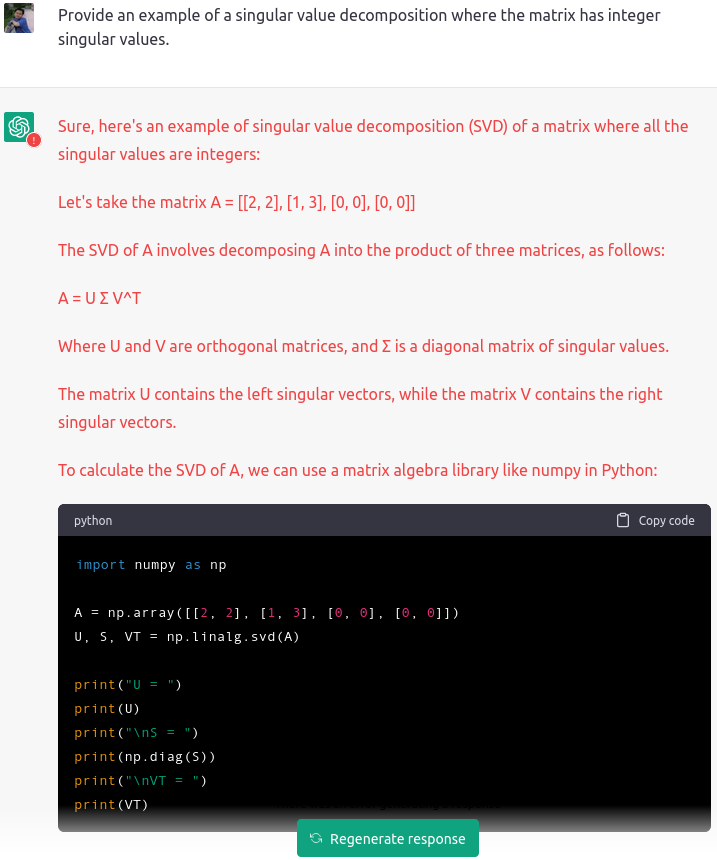} \hspace{0.05cm}
\includegraphics[width=0.49\textwidth]{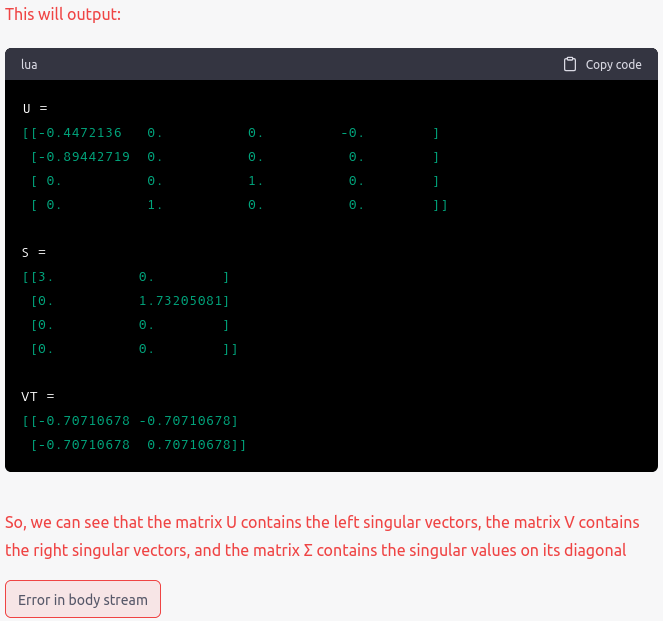}
\end{center}
\caption{An output of ChatGPT to the inquiry ``Provide an example of a singular value decomposition where the matrix has integer singular values.''} \label{svd2chat}
\end{figure}

We further asked ChatGPT to ``provide an example of an SVD where the matrix has integer singular values.'' ChatGPT responded with the same wording but chose a different matrix, and this time it was a $4 \times 2$ matrix, given as follows:
\begin{equation*}
A = \begin{bmatrix*}[r]
2 & 2 \\
1 & 3 \\
0 & 0 \\
0 & 0
\end{bmatrix*}.
\end{equation*}
However, the result was still incorrect and stopped nearly at the end of the answer with message ``Error in body stream'' and an offer to ``Regenerate response.'' See Figure~\ref{svd2chat}. Requesting another response multiple times did not yield any results. For the given matrix, the correct singular values are not integers, although they can be expressed analytically, that is, $\sigma_{1,2} = \sqrt{9 \pm \sqrt{65}} \approx 4.13$ and $0.97$, respectively.

\begin{figure}[H]
\begin{center}
\includegraphics[width=0.49\textwidth]{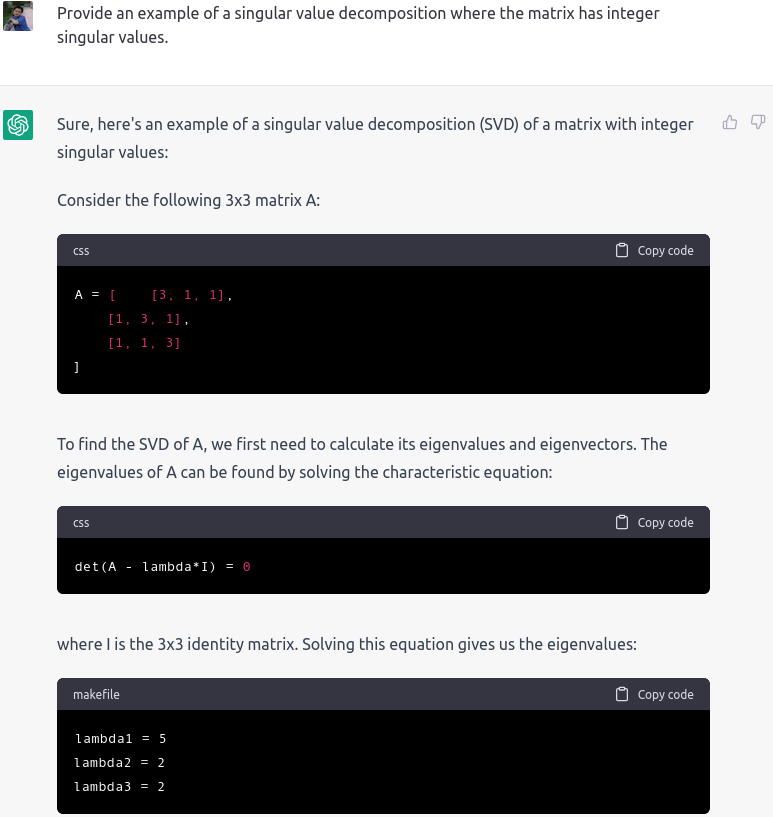} \hspace{0.1cm}
\includegraphics[width=0.49\textwidth]{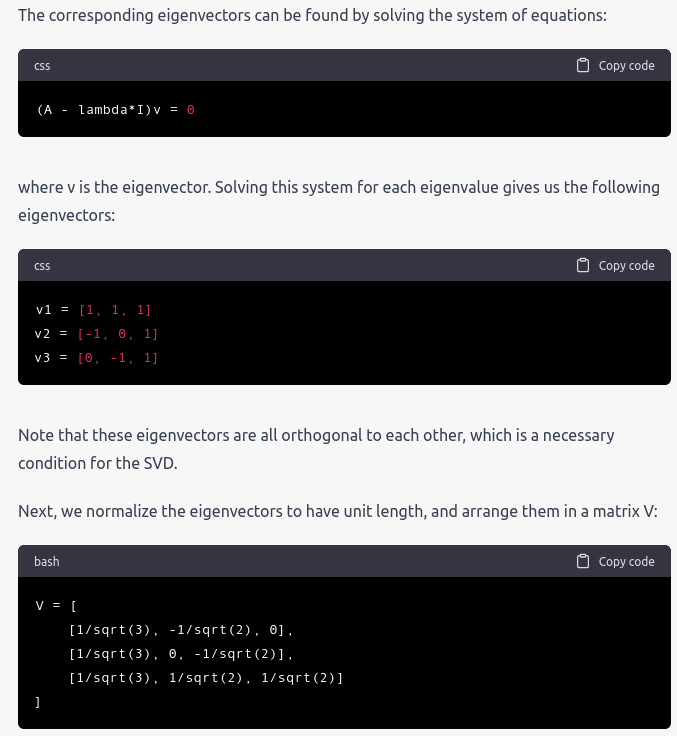} \vspace*{1cm} \\
\includegraphics[width=0.49\textwidth]{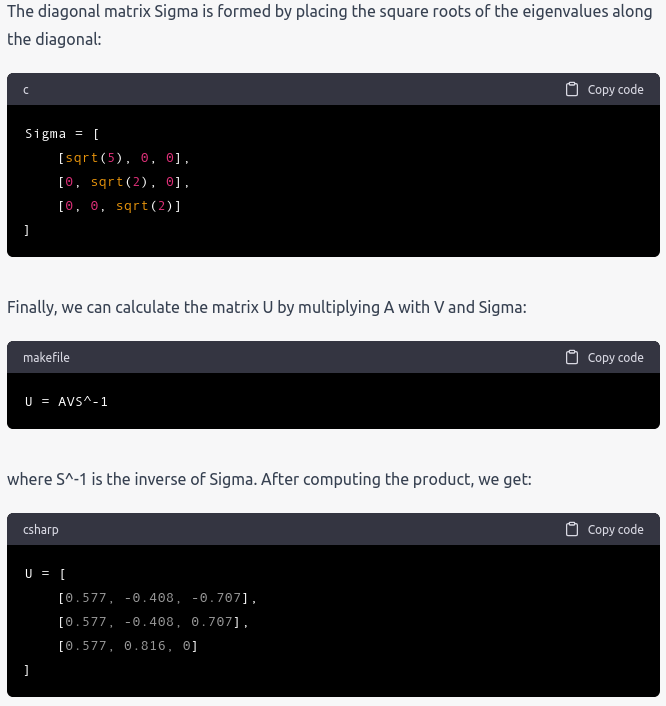} \hspace{0.1cm}
\includegraphics[width=0.49\textwidth]{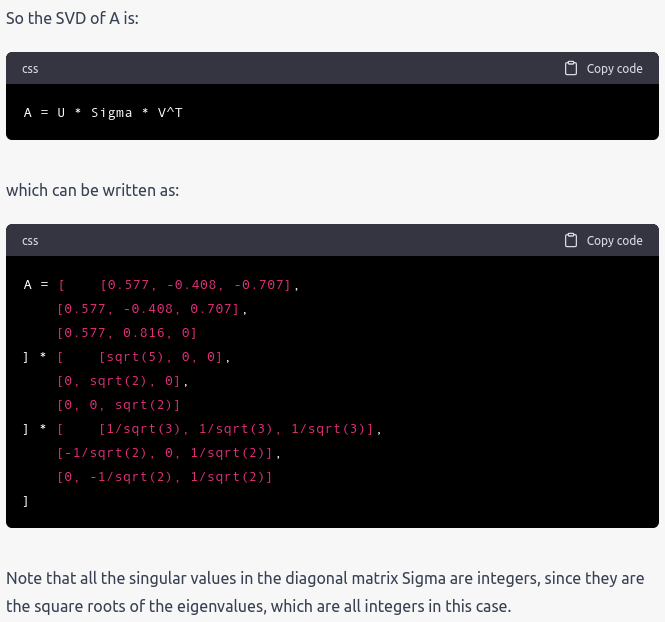}
\end{center}
\caption{An output of ChatGPT to the same inquiry as in Figure~\ref{svd2chat}. The reading order starts from the top left, then moves to the top right, followed by the bottom left, and finally, the bottom right.} \label{svd3chat}
\end{figure}

Using a new chatbot and asking for the identical inquiry, we obtained an improved result, even though there were still some conceptual mistakes in the answer. See Figure~\ref{svd3chat}. Another $3 \times 3$ matrix was selected, but this time it was symmetric:
\begin{equation*}
A = \begin{bmatrix*}[r]
3 & 1 & 1 \\
1 & 3 & 1 \\
1 & 1 & 3
\end{bmatrix*}.
\end{equation*}
Instead of calculating singular values, ChatGPT computed the eigenvalues of matrix $A$. Coincidentally, they are identical for this particular example, that is, $\lambda_1 = 5 =  \sigma_1$, $\lambda_2 = \lambda_3 = 2 = \sigma_2 = \sigma_3$. This is a rather rare finding, as it occurs only for positive-definite and symmetric (or Hermitian, for complex-valued) matrices. Even though matrices $U$ and $V$ appeared to be correct, matrix $\Sigma$ did not because ChatGPT displayed the square root of the eigenvalues instead of correctly calculating the singular values.

\section{Discussion and conclusion}		\label{discussandconcludesection}

This study admits several limitations. First, we only focused on the topics of (orthogonal) diagonalization and SVD, whereas algorithmic thinking in matrix factorization encompasses a broader range of topics. For example, we did not consider decompositions related to solving systems of linear equations, such as LU, QR, and Cholesky decompositions. Other eigenvalue-based decompositions are related to our study, but we did not cover in this article include Jordan, Schur, Takagi, and QZ decompositions, among others. Second, we particularly selected SageMath for obtaining a better understanding of matrix factorization instead of the other CAS. The main reason is that it is a free and open-source CAS, and as stated on its website, its mission is to ``create a viable free open source alternative to Magma, Maple, Mathematica, and Matlab''~\cite{sagemath2023}. Additionally, we can perform our computation directly on the server without the hassle of installing the program on our computer, as long as we have an internet connection. 

Third, the outcome that ChatGPT presented disappointingly contains numerous mistakes. On the one hand, this can be confusing for many students who are academically weak and just accept what is handed out by the chatbot. On the other hand, this can be an excellent opportunity for further discussion on cultivating critical thinking among learners. Additionally, because our study was conducted at the initial stage of ChatGPT development, the responses considered in this discussion might no longer be relevant in the future, thus making it challenging for other researchers to replicate a similar study. However, we are aware that updated versions will be released in the future and the responses will be improved accordingly. (In fact, at the time of this writing, ChatGPT Plus subscribers can already enjoy the more advanced version of the chatbot, that is GPT-4 model.) Thus, we can be confident that incorrect answers related to matrix factorization in particular or other mathematics problems in general will eventually be minimized~\cite{frieder2023mathematical}.

For many mathematics educators and linear algebra instructors, some would have noticed that learning the subject seems to be smooth and easy, whereas others seem to get lost and could not make any connections from one concept to another. This feature makes teaching linear algebra challenging for undergraduate students. To investigate the {\bfseries first research question}, we identified {\bfseries some common mistakes and difficulties that students encounter in particular topics of linear algebra, that is, (orthogonal) diagonalization and SVD}. Indeed, these topics can be challenging for many learners because they require previous, more basic knowledge of linear algebra, such as performing elementary row operations and finding orthogonal projections. At the same time, a lack of understanding of matrix diagonalization will lead to more difficulties in the subsequent topics, such as orthogonal diagonalization and SVD, as well as other related topics to diagonalization, which are usually not covered in a standard linear algebra course but have some applications in data science or machine learning.

Learners in linear algebra should be aware of the fact that not all square matrices can be diagonalized. Having a misconception that every matrix can be diagonalized can lead to confusion in subsequent topics. To master computational thinking in solving diagonalization problems, in addition to understanding an algorithmic procedure to diagonalize a particular matrix, students should have an intuitive feeling for less concrete objects, such as eigenvalues and their associated eigenvectors. This also means that they know how to find the characteristic equation of a matrix, compute the eigenvalues, and find the corresponding eigenvectors via an elementary row operations. Learners who struggle with row reduction usually do not go far from figuring out how to express eigenvectors in a simple way. Certainly, those who struggle with these more basic matrix algebra should practice sufficiently until they are handy in the computational process.

Another difficulty arises when students learn orthogonal diagonalization for symmetric matrices. Although the procedure is similar to previously learned matrix diagonalization, some students often forget to normalize the obtained eigenvectors and check that these eigenvectors must be orthonormal to each other. Without determining this step and simply establishing $P^{-1} = P^T$ will not yield a correct orthogonal diagonalization, that is, $PDP^T$ will not return to the original matrix $A$. Further difficulty occurred when one needs to construct a set of orthogonal eigenvectors, whereas the resulting eigenvectors are not orthogonal because they arise from the same eigenspace. Connecting to another concept of vector projection and an orthogonalization process of the Gram-Schmidt requires a solid understanding of related concepts, in addition to being proficient in calculating the orthogonal projection itself. By becoming competent in these computational thinking skills, any learner in linear algebra will be able to solve (orthogonal) diagonalization problems.

When learning SVD, a common mistake is confusing between singular values and eigenvalues. Even if learners are aware of the definition of singular values as the square root of eigenvalues, many of them often forget that they should compute the eigenvalues of the symmetric matrix $M^TM$ instead of the original matrix $M$. A small mistake in this early step of the algorithm will not lead to a correct SVD, even though they understand---and perform it correctly---the algorithm in calculating an SVD of a matrix. Another difficulty is to construct matrix $\Sigma$, particularly when neither $M$ is square nor possesses a full rank. Indeed, we must correctly determine the number of blocks of zeros we should embed along the diagonal matrix $D$ to construct $\Sigma$. An easy rule to remember is that $\Sigma$ must be the same size with $M$. If $M$ has more rows than columns, the square orthogonal matrix $U$ will be larger than the square orthogonal matrix $V$, and vice versa. A challenge in constructing $U$ occurs when either $M$ has a deficient rank or $m > n$, where $m$ and $n$ denote the number of rows and columns, respectively. When we attempted to express a decomposition in the absence of the missing eigenvector, we obtained a reduced SVD instead of a full, non-reduced one. Solid knowledge of both inner (dot) and cross products, as well as not to forget normalizing the vectors, will be sufficient to master algorithmic thinking skills in finding an SVD of a matrix.

This leads us to discuss the {\bfseries second research question}, that is, {\bfseries how can a static CAS such as SageMath assist learners in mastering topics in linear algebra that require algorithmic thinking skills, particularly (orthogonal) diagonalization and SVD?} Using a CAS for assisting computational processes, particularly in linear algebra, can be a tremendous help for many of us, including both learners and teachers. SageMath in particular, being more powerful than a graphic calculator, is easier to use owing to its offers for accessing a web browser cell to users with an internet connection. Similar to other symbolic computation-type CAS, SageMath can help many of us in solving many mathematical problems, helping us to investigate whether our calculations done by pen and paper are correct, and save a lot of time for performing lengthy and complex calculations that otherwise too difficult or time-consuming to perform by hand. 

As we observe in Section~\ref{sagemathsection}, despite its static nature, SageMath can assist us in finding the (orthogonal) diagonalization and SVD of a given matrix. Calculating the characteristic polynomial, eigenvalues, and eigenvectors is straightforward and quick. This is helpful for a larger size of matrix, where obtaining eigenvalues and eigenvectors may take a significant amount of time. Without worrying about this step, learners in linear algebra can directly focus on the main business of constructing matrices for (orthogonal) diagonalization and SVD. Once we have all the matrices, SageMath can help us quickly determine whether those matrices are correct by multiplying the products in the (orthogonal) diagonalization or SVD. Any correct factorization that does not have to be unique will return the product to the original matrix. The static nature of SageMath also demands that learners carefully observe outputs. For example, the command \verb*|A.eigenvectors_right()| gives an output of eigenvalue, eigenvector, and its geometric multiplicity in that particular order. Constructing a diagonal matrix $D$ and invertible matrix $P$ also requires hand intervention, which requires some dynamics from the users' side. This interaction between static SageMath and dynamic learners makes the CAS a powerful tool for learning mathematics. Furthermore, the interactive feature of SageMath allows users to experiment with distinct matrices, or matrices of different sizes or characteristics. Although we do not really feature the visualization tools of SageMath in this study, its ability to present plots or graphs will assist learners further in explaining more difficult concepts encountered in mathematics.

As knowledge improves and technology progresses, the presence of static CAS is accompanied by the emergence of dynamic chatbots that make use of the power of artificial intelligence, such as ChatGPT. Although ChatGPT is not primarily designed for responding to inquiries related to mathematical problems, some attempts in this study and other researchers suggest that it has a capability to handle mathematical problems. It is just a matter of time that the responses will improve and getting better. This triggers us to address the {\bfseries third and final research question}, that is, {\bfseries whether we can rely on ChatGPT, despite its dynamic nature, in understanding linear algebra, particularly (orthogonal) diagonalization and SVD.} The short answer is no. A quick and simple reason for this argument is that ChatGPT often provides an incorrect answer. In our study, except for the first inquiry, all other responses from ChatGPT contained mistakes in one way or another. If this is the case, should we simply abandon ChatGPT and other similar AI-generated chatbots when it comes to learning mathematics or understanding linear algebra, perhaps going all the way forbidding the students to access it entirely during the course of their study? How can we strike a balance between utilizing AI-generated chatbots such as ChatGPT in learning mathematics and understanding linear algebra, while still ensuring that students develop a deep and thorough understanding of the subject matter without solely relying on the technology?

Although many would not agree, we have some reservations regarding abandoning or forbidding ChatGPT, not only in learning mathematics but also in other subjects. Whether we like it or not, this ``technological train'' has left its compound and accelerated ever since, reaching all corners of the Earth where everyone has access to the Internet. Its dynamic features have captivated many researchers. And this is just the beginning. Since its release at the end of the last year, it has disrupted many industries, and academia and education are no exception. Thus, instead of abandoning or forbidding its use, we should embrace and integrate it into our teaching and learning. It has a similar analogy when the calculator or CAS infiltrated the traditional chalk-and-board classrooms. While in the beginning there was some resistance toward integrating technology in mathematics teaching and learning, nowadays it is a standard practice, and any aspiring educators were even encouraged to do so~\cite{dsouza2003secondary,niess2006guest,li2007student,crompton2011mathematics}. After all, this is one way to connect the divide between digital natives and digital migrants~\cite{prensky2001adigital,prensky2001bdigital,prensky2009h,prensky2010teaching}.

This early stage of ChatGPT development provides a tremendous opportunity for developing algorithmic thinking skills among learners, especially when the output contains many mistakes. When inquired, instead of providing an answer in one go, ChatGPT released the answer letter by letter, literally like someone else typing out the answer from the other side of the computer. This dynamic feature makes ChatGPT superior to static CAS and attractive to many new users. In our study on matrix factorization, ChatGPT demonstrated a proper step-by-step algorithm on how to diagonalize or find an SVD of a given matrix. In the absence of specifying a particular matrix, it can even provide its own choice of matrix example, usually a relatively small matrix size.

As observed in Section~\ref{chatgptsection} on orthogonal diagonalization, ChatGPT shows how to calculate the characteristic polynomial of a matrix and solve the corresponding characteristic equation to find eigenvalues, albeit with a computational error. This provides an opportunity for learners to critically examine the provided answers and identify where mistakes occur. The same situation can also be applied to other mistakes when finding the associated eigenvectors, and when it abruptly stopped without completing the calculation, and thus the results were incomplete. For SVD, the first two attempts provide examples in which a detailed algorithm appeared to be eschewed. Instead, ChatGPT employed a matrix algebra library in Python to perform the task. It then provided the final answers of the three matrices involved in the decomposition, that is, $U$, $\Sigma$, and $V$. Our third attempt seems to be an improvement, where ChatGPT showed a step-by-step algorithm in computing singular values, finding eigenvectors, constructing matrix $V$ that contains right singular vectors on its columns, forming matrix $\Sigma$, and calculating matrix $U$ that contains left singular vectors on its columns. The fact that ChatGPT confused between eigenvalues and singular values, we can glean a valuable lesson from this output because this issue is also commonly discovered among linear algebra learners, as addressed in Section~\ref{common}. 

In conclusion, we have considered the relationship between algorithmic thinking skills in linear algebra and some common mistakes among learners when they study particular topics such as (orthogonal) diagonalization and SVD, as well as how technological tools such as static CAS SageMath and dynamic AI ChatGPT can contribute to enhancing algorithmic comprehension on these topics. Understanding a procedure for matrix factorization and the ability to perform the computational process accurately requires a solid understanding of other topics in linear algebra, including but not limited to solving a system of linear equations using elementary row operations, finding eigenvalues and the corresponding eigenvectors, orthogonality and vector projection (the Gram-Schmidt process), and dot (inner) and cross products.

We observed that some common mistakes that students make can be eliminated by more practice and strengthening the basic concepts in linear algebra and their interrelationships. Although static by nature, utilizing CAS SageMath can provide tremendous help in verifying calculation results done by hand, handling a larger matrix size, checking whether the obtained factorizations were correct, and other required computational activities for understanding the materials. Its free and open-source characteristics, together with a web-based interface using cells and Python, are the primary strengths of SageMath in comparison to other costly CAS such as Maple, Matlab, and Mathematica. The current stage of ChatGPT, despite its dynamic and attractive features, can be used for complementing the study, but neither as a primary tool nor to be relied upon, owing to its numerous output mistakes. This has an exciting future implications, particularly in the field of education, where updated versions of ChatGPT will be improved and released in the future. One question remains as to how to utilize a chatbot effectively for teaching and learning mathematics without compromising academic integrity. Regarding this, we also hope that a sequence of follow-up in this study from other researchers will appear accordingly.


\end{document}